\documentclass[11pt,reqno]{amsart}
\usepackage{amssymb, hyperref}

\numberwithin{equation}{section}
\newtheorem{theorem}{Theorem}[section]
\newtheorem{proposition}[theorem]{Proposition}
\newtheorem{lemma}[theorem]{Lemma}

\newtheorem{remark}[theorem]{Remark}

\newcommand{\R}{\ensuremath{\mathbb{R}}}

\newcommand{\px}{\partial_x}
\newcommand{\py}{\partial_y}

\newcommand{\Nn}{N_{min}}
\newcommand{\Nd}{N_{med}}
\newcommand{\Nx}{N_{max}}
\newcommand{\Ln}{L_{min}}
\newcommand{\Ld}{L_{med}}
\newcommand{\Lx}{L_{max}}

\begin{document}
\title[fifth order modified KdV equation]
  { Well-posedness and ill-posedness of the fifth order modified KdV equation}
\author[Soonsik Kwon ]
{Soonsik Kwon}

\address{Soonsik Kwon \hfill\break
Department of Mathematics, University of California, Los Angeles, CA
90095-1555, USA} \email{rhex2@math.ucla.edu}

\begin{abstract}
We consider the initial value problem of the fifth order modified KdV equation on the Sobolev spaces.
\begin{equation*}
\begin{cases}
\partial_t u - \px^5u + c_1\px^3(u^3) + c_2u\px u\px^2 u + c_3uu\px^3 u =0\\
u(x,0)= u_0(x)
\end{cases}
\end{equation*}
where $ u:\R\times\R \rightarrow \R $ and $c_j$'s are real. We show the local well-posedness in $H^s(\R)$ for $s\geq \frac{3}{4}$ via the contraction principle on $X^{s,b}$ space. Also, we show that the solution map from data to the solutions fails to be uniformly continuous below $H^{3/4}(\R)$. The counter example is obtained by approximating the fifth order mKdV equation by the cubic NLS equation.
\end{abstract}

\thanks{}
\thanks{} 

\maketitle

\section{Introduction}

The KdV equation and the modified KdV(mKdV) equation are completely integrable in the sense that there are Lax pair formulations. Being completely integrable, the KdV and the mKdV equations enjoy infinite number of conservation laws. Each of these is an Hamiltonian of the flow which commute the KdV flow (resp. the mKdV flow). This generates an infinite collection of commuting nonlinear equations of order $2j+1,\,\,(j\in \mathbb{N})$, which is known as the KdV hierarchy (resp. the mKdV hierarchy). In this note, we consider the second equation from the modified KdV hierarchy:
\begin{equation}\label{fifthmkdv}
\partial_t u - \px^5 u -30u^4\px u + 10 u^2\px^3 u + 10 (\px u)^3 + 40 u\px u\px^2 u = 0.
\end{equation}
Using the theory of the complete integrability, one can show that for any Schwartz initial data, the solution to any equation in the KdV hierarchy (resp. the mKdV hierarchy) exists globally in time. However, the well-posedness theory for low regularity initial data requires the theory of dispersive PDEs. And changing coefficients in the nonlinear terms may break the integrable structure. In this case, we can no longer rely on the theory of complete integrability. The purpose of this paper is to study the low regularity well-posedness and ill-posedness.\\
We consider the following fifth order mKdV equation, which generalizes \eqref{fifthmkdv} \footnote{For omitting $u^4\px u$, See Remark~\ref{low order term}}.
\begin{equation}\label{fifthm}
\begin{cases}
\partial_t u - \px^5u + c_1\px^3(u^3) + c_2u\px u\px^2 u + c_3uu\px^3 u =0\\
u(x,0)= u_0(x)
\end{cases}
\end{equation}
where $ u:\R\times\R \rightarrow \R $ and $c_j$'s are real. \\
\\
We show the local well-posedness result and the ill-posedness result.
First, we state the local well-posedness theorem.
\begin{theorem}\label{lwp}
Let $ s \geq \frac 34$ and $u_0 \in H^s(\R)$. Then there exists $T=T(\|u_0\|_{H^s(\R)})$ such that the initial value problem \eqref{fifthm} has a unique solution $u(t,x) \in C([0,T];H^s(\R))$. Moreover, the solution map from data to the solutions is real-analytic.
\end{theorem}
Previously, Kenig, Ponce, and Vega \cite{KPV94} studied the local well-posedness of the odd order dispersive equations:
$$ \partial_t u + \px^{2j+1}u + P(u,\px u,\cdots ,\px^{2j}u) =0 $$
where $P$ is a polynomial having no constant and linear terms. They proved the local well-posedness for the initial data in the weighted Sobolev space, i.e. $$ u_0 \in H^s(\R) \cap L^2(|x|^mdx) $$ for some $s,m \geq 0$. Their method was the iteration using the local smoothing estimate and the maximal function estimate. Inspecting their proof for the equation \eqref{fifthm}, one can observe that the local well-posedness holds true for $s>\frac{9}{4}$ and $ m=0$. In other words, the local well-posedness is established for the Sobolev space \emph{without the decaying weight}. Thus, our result can be viewed as an improvement of theirs. Our proof is also via the contraction principle. A natural choice of the iteration space is the Bourgain space, also known as the $X^{s,b}$ space. Assuming the standard argument of the iteration on the $X^{s,b}$ space, the main step is to show the following nonlinear estimate:
$$  \|T(u,v,w)\|_{X^{s,b-1}} \lesssim \|u\|_{X^{s,b}}\|v\|_{X^{s,b}}\|u\|_{X^{s,b}}. $$
where $ T(u,v,w) = c_1\px^3(uvw) + c_2u\px v\px^2 w + c_3uv\px^3 w $. This is performed by the dyadic method of Tao. In \cite{Tao2001}, Tao studied multilinear estimates for $X^{s,b}$ space systematically and showed the analogous trilinear estimate for the mKdV equation. This reproves the local well-posedness for $s\geq 1/4$, which originally showed by Kenig, Ponce and Vega \cite{KPV93} by the local smoothing estimate. Thus, in the mKdV equation the $X^{s,b}$ estimate has the same strength as the classical local smoothing method, while in the fifth order mKdV \eqref{fifthm} the $X^{s,b}$ estimate improves the preceding one. \\
\indent In \cite{CCT2003} Christ, Colliander and Tao showed the solution map of the mKdV equation fails to be uniformly continuous for $s <1/4$. This implies $1/4$ is the minimal regularity threshold for which the well-posedness problem can be solved via an iteration methods. Our next theorem is the analogue of this for the equation \eqref{fifthm}.
\begin{theorem}\label{illposed}
Let $  -\frac {7}{24} < s < \frac 34$. The solution map of the initial value problem \eqref{fifthm} fails to be uniformly continuous. More precisely, for $ 0 < \delta \ll \epsilon \ll 1$ and $T>0$ arbitrary, there are two solutions $ u,v$ to \eqref{fifthm} such that
\begin{gather}
\label{illposed1}\|u(0)\|_{H^s_x}, \|v(0)\|_{H^s_x} \lesssim \epsilon\\
\label{illposed2}\|u(0)-v(0)\|_{H^s_x} \lesssim \delta \\
\label{illposed3}\sup_{0\leq t\leq T} \|u(t) -v(t)\|_{H^s_x} \gtrsim \epsilon.
\end{gather}
\end{theorem}
The method used here is very similar to theirs in \cite{CCT2003}. We approximate the fifth order mKdV equation by the cubic NLS equation,
\begin{equation}\label{cNLS}
i\partial_t u - \px^2 u + |u|^2u = 0,
\end{equation}
at $(N,N^5)$ in the frequency space. \\ Let $u(t,x)$ be the linear solution to $ (\partial_t -\px^5)u =0 $ with $u(0)=u_0$. Setting
$$ \xi := N+ \frac{\xi'}{\sqrt{10}N^{3/2}}, $$
$ \tau=\xi^5 $ leads $ \tau = N^5 + \sqrt{\frac 52}N^{5/2}\xi' + \tau' $ where
$$ \tau' = \xi'^2 + \frac{\xi'^3}{\sqrt{10}N^{5/2}} + \frac{\xi'^4}{20N^5} +\frac{\xi'^5}{(10N^3)^{5/2}}. $$
\begin{align*}
u(t,x) &= \int e^{it\tau + ix\xi} \widehat{u_0}(\xi) d\tau d\xi \\
 &= \int e^{it(N^5 + \sqrt{\frac 52}N^{5/2}\xi' + \tau') + ix(N + \frac{\xi'}{\sqrt{10}N^{3/2}})} \widehat{u_0}(\xi) d\tau d\xi\\
 &= e^{iN^5 +iNx}\int e^{i\tau' t + i\xi'( \frac{x}{\sqrt{10}N^{3/2}}+ \sqrt{\frac{5}{2}}N^{5/2}t)}\widehat{u_0}(N+ \frac{\xi'}{\sqrt{10}N^{3/2}}) d\tau' d\xi'
 \end{align*}
Since $\tau' \approx \xi'^2 $ for $|\xi'| \ll N $,
$$u(t,x)  \approx e^{iN^5t +iNx} v( t, \frac{x}{\sqrt{10}N^{3/2}} + \sqrt{\frac{5}{2}}N^{5/2}t) $$
where $v(t,x) $ is a solution to the linear Schr\"odinger equation $i\partial_tv -\px^2v =0 $.\\
At the presence of the nonlinear term, one need the factor $ \frac{c}{N^{3/2}}$ and the real part projection Re. Then it is approximated to the cubic NLS equation \eqref{cNLS}. The detail follows in Section~\ref{sec4}. \\
\indent On the other hand, the solutions to the fifth order KdV equation,
$$ \partial_tu -\px^5u + c_1\px u\px^2u + c_2 u\px^3 u = 0, $$
is known to have genuine nonlinear dynamics for all $s>0$. In \cite{Kwon} the author showed the solution map fails to be uniformly continuous in $H^s(\R)$ for $s>0$. Thus, for this equation the local well-posedness problem is solved by other than the iteration method. In \cite{Kwon} the local well-posedness in $H^s(\R)$ for $s>\frac{5}{2}$ is established via the compactness method. \\
\subsection*{Notation}
We use $X\lesssim Y$ when $X \leq CY $ for some $C$. We use $ X \sim  Y $ when $ X \lesssim Y $ and $ Y\lesssim X$.  Moreover, we use $ X \lesssim_s Y $ if the implicit constant depends on $s$, $C=C(s) $. \\
We use Japanese bracket notation $\langle\xi\rangle := \sqrt{1+\xi^2} $.
We denote the space time Fourier transform by $\widetilde{u}(\tau,\xi) $ of $u(t,x)$
$$ \widetilde{u}(\tau,\xi) = \int e^{-it\tau - ix\xi} u(t,x) dtdx, $$
while the space Fourier transform by $\widehat{u}(t,\xi)$ of $u(t,x)$
$$ \widehat{u}(t,\xi)= \int e^{-ix\xi} u(t,x) dx.$$

\subsection*{Acknowledgement}
The author would like to appreciate his advisor Terence Tao for many helpful conversations and encouragement.

\section{Local well-posedness of the fifth order modified KdV}

In this section, we prove the local well-posedness of the initial value problem \eqref{fifthm}. Our proof is via the contraction principle on the Bourgain space. We first recall some standard facts and notations. For a Schwartz function $u_0(x) $, we denote the linear solution $u(t,x)$ to the equation $ \partial_t u - \partial_x^5 u = 0 $ by
$$ u(t,x) =: e^{t\px^5}u_0(x) =  c \iint e^{it\xi^5}e^{i(x-y)\xi}u_0(y) dyd\xi. $$
Using this notation we have the Duhamel formula for the solution to the inhomogeneous linear equation $ \partial_t u - \partial_x^5 u + F = 0 $
$$ u(t,x) = e^{t\px^5}u_0(x) - \int_0^t e^{(t-t')\px^5}F(t',x) dt'.  $$
We denote the Bourgain space by $X^{s,b}_{\tau=\xi^5}(\R\times\R)$, or abbreviated $X^{s,b}$. The $X^{s,b}$ space is defined to be the closure of the Schwartz functions $\mathcal{S}(\R\times\R) $ under the norm
$$ \|u\|_{X^{s,b}_{\tau=\xi^5}(\R\times\R)} := \|\langle \xi\rangle^s\langle\tau-\xi^5\rangle^b\widetilde{u}(\tau,\xi)\|_{L^2_{\tau,\xi}(\R\times\R)}. $$
The $X^{s,b}$ space is continuously embedded in $C^0_tH^s_x$.
\begin{lemma}
Let $b>1/2$ and $s\in \R$. Then for any $u\in X^{s,b}_{\tau=\xi^5}(\R\times\R)$, we have
$$ \|u\|_{C^0_tH^s_x(\R\times\R)} \lesssim_b \|u\|_{X^{s,b}_{\tau=\xi^5}(\R\times\R)}. $$
\end{lemma}
For the proof see \cite{Taobook}. \\
Let $\eta(t)$ be a compactly supported smooth time cut-off function (i.e. $\eta \in C^\infty_0(\R)$ with $\eta(t)=1$ on $[0,1]$). There is a standard $X^{s,b}$ energy estimate for time cut-off solutions.
\begin{lemma}
Let $b>1/2$ and $s\in \R$ and let $u\in C^\infty_t\mathcal{S}_x(\R\times\R)$ solves the inhomogeneous linear fifth order KdV equation $\partial_tu-\px^5u =F$. Then we have
\begin{equation}\label{Xlinear}
\|\eta(t)u\|_{X^{s,b}_{\tau=\xi^5}(\R\times\R)} \lesssim_{\eta,b} \|u(0)\|_{H^s_x} + \|\eta(t)F\|_{X^{s,b-1}_{\tau=\xi^5}(\R\times\R)}.
\end{equation}
\end{lemma}
For the proof of this Lemma, See \cite{KPV96}, \cite{Taobook}.\\
Next, we state the nonlinear estimate.
\begin{proposition}\label{trilinear}
\,Let $s\geq \frac 34$. For all $u,v,w$ on $\R\times \R$ and $\frac 12 <b\leq \frac 12 +\epsilon $ for some $\epsilon$, we have
\begin{equation}\label{trilinear1}
\begin{split}
\|\px^3(uvw)\|_{X^{s,b-1}} &+ \|uv\partial_x^3w\|_{X^{s,b-1}}+\|u\px^2v\px w\|_{X^{s,b-1}} +\|\px u\px v\px w\|_{X^{s,b-1}} \\
&\lesssim \|u\|_{X^{s,b}}\|v\|_{X^{s,b}}\|w\|_{X^{s,b}}.
\end{split}
\end{equation}
\end{proposition}

Combining the preceding estimates \eqref{Xlinear}, \eqref{trilinear1} one can easily verify that the operator
$$ \Phi(u)(t,x) := \eta(t)e^{t\px^5}u_0(x) - \eta(t)\int^t_0 e^{(t-t')\px^5}F(u)(t',x) dt' $$
is a contraction on a ball of $X^{s,b}$ space
$$ \mathcal{B}= \{u \in X^{s,b} : \|u\|_{X^{s,b}} \leq 2\delta\} $$
for a sufficiently small $ \delta >0 $ and $ \|u_0\|_{H^s_x} <\delta $, where $F(u) = c_1\px^3(u^3) + c_2u\px u\px^2 u + c_3uu\px^3 u$. This proves the local well-posedness for small data. Then a standard scaling argument easily leads the local well-posedness for large data. Once the local well-posedness is proved via the contraction principle, we also obtain that the solution map is Lipschitz continuous, and furthermore if the nonlinear term is algebraic (a polynomial of u and its derivatives), then the solution map is real-analytic. Hence, it remains to show the trilinear estimate \eqref{trilinear1} for the proof of Theorem~\ref{lwp}.

\section{Trilinear estimate}
In this section, we show the trilinear estimate \eqref{trilinear1}. We closely follow the method developed by Tao \cite{Tao2001} in the context of modified KdV equation. Writing the trilinear inequality in the dual form and we view it as a composition of two bilinear operators based on $L^2$ norm. Then we reduce to two bilinear estimates. First, we recall notations and general frame work of Tao's $[k;Z]$-multiplier method. For the details we refer to \cite{Tao2001}.
\subsection*{Notation and block estimates} We define $[k,\R]$-multiplier norm of Tao \cite{Tao2001} first. Let $Z$ be an abelian additive group with an invariant measure $d\xi$ (for instance $\R^n, \mathbb{T}^n$). For any integer $k\geq2$. let $\Gamma_k(Z)$ denote the hyperplane
$$ \Gamma_k(Z) := \{(\xi_1,\cdots,\xi_k) \in \R^k : \xi_1+\cdots+\xi_k =0 \}. $$
A $[k,Z]$-multiplier is defined to be any function $m: \Gamma_k(Z) \rightarrow \mathbb{C}$. Then we define the multiplier norm $\|m\|_{[k,Z]} $ to be the best constant so that the inequality
$$ \big| \int_{\Gamma_k(Z)} m(\xi)\prod_{j=1}^{k}f_j(\xi_j)\big| \leq C \prod_{j=1}^k\|f_j\|_{L^2},
 $$ holds for all functions $f_j$ on $Z$. \\
Any capitalized variables such as $N_j,L_j$ and $H$ are presumed to be dyadic. For $N_1,N_2,N_3>0$, we denote the quantities by $\Nn,\Nd,\Nx$ in their order and similarly for $L_1,L_2,L_3$. We adopt the following summation convention:
$$ \sum_{\Nx\sim\Nd\sim N} := \sum_{\substack{ N_1,N_2,N_3>0\\ \Nx\sim\Nd\sim N}} $$
$$ \sum_{\Lx\sim H} := \sum_{\substack{L_1,L_2,L_3\gtrsim 1\\ \Lx\sim H}}. $$
For given $\tau_j, \xi_j$ with $\xi_1 +\xi_2 +\xi_3 =0$ and $\tau_1+\tau_2+\tau_3=0$, we denote the modulation $$\tau_j-\xi_j^5 =:\lambda_j$$ and the resonance function
$$ h(\xi):= \xi_1^5 +\xi_2^5 +\xi_3^5 = -\lambda_1 -\lambda_2 -\lambda_3. $$
By a dyadic decomposition of the variables $\xi_j, \lambda_j$ and $h(\xi)$ $X^{s,b}$, a bilinear estimate
$$ \|B(u,v)\|_{X^{s_3,b_3}} \lesssim \|u\|_{X^{s_1,b_1}}\|v\|_{X^{s_2,b_2}} $$
is reduced to
$$ \left\| \sum_{\Nx\gtrsim 1}\sum_{H}\sum_{L_1,L_2,L_3\gtrsim 1} \frac{\widetilde{m}(N_1,N_2)\langle N_1\rangle^{-s_1}\langle N_2\rangle^{-s_2}\langle N_3\rangle^{s_3}}{L_1^{b_1}L_2^{b_2}L_3^{-b_3}}X_{N_1,N_2,N_3;H;L_1,L_2,L_3}  \right\|_{[3,\R\times\R]}  \lesssim 1.  $$
Here, $X_{N_1,N_2,N_3;H;L_1,L_2,L_3}$ is the multiplier
$$ X_{N_1,N_2,N_3;H;L_1,L_2,L_3}(\xi,\tau) := \chi_{|h(\xi)|\sim H}\prod_{j=1}^3 \chi_{|\xi_j|\sim N_j}\chi_{|\lambda_j|\sim L_j}  $$
and
$$ \widetilde{m}(N_1,N_2) := \sup_{|\xi_j|\sim N_j, j=1,2} m(\xi_1,\xi_2)$$
where $ m(\xi_1,\xi_2) $ is a multiplier of the bilinear operator $B(\cdot,\cdot)$.
This leads us to consider
\begin{equation}\label{block}
\|X_{N_1,N_2,N_3;H;L_1,L_2,L_3}\|_{[3,\R \times\R]},
\end{equation}
which vanishes unless
\begin{gather}
\label{Nmed-Nmax} \Nd \sim \Nx\\
\label{Lmax} \Lx \sim \max(H, \Ld)
\end{gather}
Moreover, we have the resonance relation: if $ \Nx \sim \Nd \gtrsim 1$, then
\begin{equation}\label{resonance}
H \sim \Nx^4\Nn
\end{equation}
Now we state the dyadic block estimate.

\begin{lemma}\label{block estimate} Let
$H,\,N_1,\,N_2,\,N_3,\,L_1,\,L_2,\,L_3>0$ obey \eqref{Nmed-Nmax}, \eqref{resonance}, \eqref{Lmax}.

(a)((++)Coherence) If $N_{max}\sim N_{min}$ and $L_{max}\sim
H$, then we have
\begin{eqnarray}\label{estimate1}
   \eqref{block} \lesssim L_{min}^{1/2}N_{max}^{-2}L_{med}^{1/2}.
\end{eqnarray}

(b)((+-)Coherence) If $N_2\sim N_3\gg N_1$ and $H\sim
L_1\gtrsim L_2,\,L_3$, then
\begin{eqnarray}\label{estimate2}
 \eqref{block}\lesssim
L_{min}^{1/2}N_{max}^{-2}\min(H,\,\frac{N_{max}}{N_{min}}L_{med})^{1/2}.
\end{eqnarray}
Similarly for permutations.

(c) In all other cases, we have
\begin{eqnarray}\label{estimate3}
\eqref{block}\lesssim
L_{min}^{1/2}N_{max}^{-2}\min(H,\,L_{med})^{1/2}.
\end{eqnarray}
\end{lemma}
Lemma~\ref{block estimate} is obtained in a similar way to Tao's (\cite{Tao2001}, Proposition 6.1) in the context of the KdV equation. For the fifth order equation, it is first shown by Chen, Li, Miao and Wu \cite{Chen-et-al}. See \cite{Chen-et-al} for the proof.

\subsection*{Bilinear estimates}
Using Lemma~\ref{block estimate} we show three bilinear estimates to which the trilinear estimate is reduced.
\begin{proposition}\label{bilinear}
For Schwartz functions $u, v$ on $\R\times \R$ and $ 0<\epsilon \ll 1$, we have
\begin{align}
\label{bilinear1}
\|uv\|_{L^2(\R\times\R)} &\lesssim \|u\|_{X^{-3/2,1/2-\epsilon}_{\tau=\xi^5}}\|v\|_{X^{3/4,1/2+\epsilon}_{\tau=\xi^5}}, \\
\label{bilinear2}
\|uv\|_{L^2(\R\times\R)} &\lesssim \|u\|_{X^{-3/4,1/2-\epsilon}_{\tau=\xi^5}}\|v\|_{X^{0,1/2+\epsilon}_{\tau=\xi^5}}, \\
\label{bilinear3}
\|uv\|_{L^2(\R\times\R)} &\lesssim \|u\|_{X^{-1/4,1/2-\epsilon}_{\tau=\xi^5}}\|v\|_{X^{-1/2,1/2+\epsilon}_{\tau=\xi^5}}.
\end{align}

\end{proposition}

\begin{proof}
We prove \eqref{bilinear1} first. Rewriting \eqref{bilinear1} by duality, Plancherel's theorem and dyadic decomposition and using the translation invariance of the $[k;Z]$-multiplier (may assume $L_1,L_2,L_3 \gtrsim 1$ and $\max(N_1,N_2,N_3) \gtrsim 1$) and Schur's test (\cite{Tao2001}, Lemma 3.11), it suffices to show
\begin{align}\label{H-Lmax}
\sum_{N\sim \Nx \sim \Nd} \sum_{\substack{L_1,L_2,L_3\geq 1,\\ H\sim \Lx}} \\
\frac{\langle N_2\rangle^{3/2}}{\langle N_1 \rangle^{3/4}L_1^{1/2+\epsilon}L_2^{1/2-\epsilon}}
 &\|X_{N_1,N_2,N_3;\Lx;L_1,L_2,L_3}\|_{[3,\R \times\R]} \lesssim 1  \nonumber
\end{align}
and
\begin{align}\label{H<Lmax}
\sum_{N\sim \Nx \sim \Nd} \sum_{\substack{\Lx \sim \Ld,\\ H\ll \Lx}} \\
\frac{\langle N_2\rangle^{3/2}}{\langle N_1 \rangle^{3/4}L_1^{1/2+\epsilon}L_2^{1/2-\epsilon}}
&\|X_{N_1,N_2,N_3;\Lx;L_1,L_2,L_3}\|_{[3,\R \times\R]} \lesssim 1 \nonumber
\end{align}
for all $N \gtrsim 1$.
Fix $N$. We prove \eqref{H<Lmax} first. From \eqref{estimate3} it reduces to show
\begin{equation*}
\sum_{N\sim \Nx \sim \Nd} \sum_{\Lx \sim \Ld \gtrsim N^4\Nn} \frac{\langle N_2\rangle^{3/2}}{\langle N_1 \rangle^{3/4}L_1^{1/2+\epsilon}L_2^{1/2-\epsilon}}
\Ln^{1/2}N^{-2}N^{2}\Nn^{1/2}  \lesssim 1
\end{equation*}
Estimating
\begin{align*}
\frac{\langle N_2 \rangle^{3/2}}{\langle N_1 \rangle^{3/4}} &\lesssim \frac{N^{3/2}}{\langle \Nn \rangle^{3/4}} \\
L_1^{1/2+\epsilon}L_2^{1/2-\epsilon} &\gtrsim \Ln^{1/2+\epsilon}\Ld^{\epsilon}(N^4\Nn)^{1/2-2\epsilon}
\end{align*}
and then performing the $L$ summations, we reduce to
\begin{equation*}
\sum_{N\sim \Nx \sim \Nd}\frac{\langle N\rangle^{3/2}\Nn^{1/2}}{\langle \Nn \rangle^{3/4}(N^4\Nn)^{1/2-\epsilon}}  \lesssim 1,
\end{equation*}
which is true with about $N^{-1/2}$ to spare.\\
Now, we show the case \eqref{H-Lmax}. In this case we have $ \Lx \sim \Nx^4\Nn $. We first show when \eqref{estimate1} (i.e. (++)coherence) holds. From \eqref{estimate1} we have $\Nx \sim \Nd \sim \Nn$ and $ \eqref{block} \lesssim \Ln^{1/2}\Nx^{-2}\Ld^{1/2} $, so we reduce to
\begin{equation*}
\sum_{\Lx\sim N^5} \frac{N^{3/2}}{N^{3/4}L_1^{1/2-\epsilon}L_2^{1/2+\epsilon}} \Ln^{1/2}N^{-2}\Ld^{1/2} \lesssim 1.
\end{equation*}
Estimating
$$ L_1^{1/2+\epsilon}L_2^{1/2-\epsilon} \geq \Ln^{1/2+\epsilon}\Ld^{1/2-\epsilon} $$ and then performing the $L$ summations we reduce to
$$  \frac{N^{3/2}}{N^{3/4}N^{5\epsilon}} N^{-2} \lesssim 1, $$
which is true.\\
\indent Now we deal with (+-)coherence case (i.e. when \eqref{estimate2} holds true). Since we don't have the symmetry on indices, we need to consider the following three cases:
\begin{align*}
N \sim N_1 \sim N_2 \gg N_3; \quad H\sim L_3 \gtrsim L_1,L_2 \\
N \sim N_2 \sim N_3 \gg N_1; \quad H\sim L_1 \gtrsim L_2,L_3 \\
N \sim N_1 \sim N_3 \gg N_2; \quad H\sim L_2 \gtrsim L_1,L_3 \\
\end{align*}
In the first case we reduce by \eqref{estimate2} to
\begin{equation*}
\sum_{N_3 \ll N, L_1,L_2 \lesssim N^4N_3} \frac{N^{3/2}}{N^{3/4}L_1^{1/2+\epsilon}L_2^{1/2-\epsilon}}\Ln^{1/2}N^{-2}\min(N^4N_3,\frac{N}{N_3}\Ld)^{1/2} \lesssim 1.
\end{equation*}
Performing the $N_3$ summation we reduce to
\begin{equation*}
\sum_{1 \leq L_1,L_2 \lesssim N^5} \frac{N^{3/2}}{N^{3/4}L_1^{1/2+\epsilon}L_2^{1/2-\epsilon}}\Ln^{1/2}N^{-2}N^{5/4}\Ld^{1/4} \lesssim 1
\end{equation*}
which is easily verified. \\
To symmetrize the second and third case we replace $L_1^{1/2+\epsilon} $ by $L_1^{1/2-\epsilon}$. It suffices to show the second case.
Using $ \min(H,\frac{N}{\Nn}\Ld) \leq H \sim N^4N_1 $ we reduce to
$$ \sum_{N_1\leq N}\sum_{L_2,L_3\leq N^4N_1} \frac{N^{3/2}N_1^{1/2}}{\langle N_1 \rangle^{3/4}(N^4N_1)^{1/2+\epsilon}L_2^{1/2-\epsilon}} \lesssim 1 $$
We may assume $ N_1 \geq N^{-4} $ since the inner sum vanishes otherwise. Performing the $L$ summations we reduce to
$$ \sum_{N^{-4} \leq N_1 \leq N} N^{3/2-2+4\epsilon}\frac{N_1^\epsilon}{\langle N_1 \rangle^{3/4}}(N^4N_1)^\epsilon \lesssim 1 $$
which is true with about $N^{-1/2}$ to spare.
Finally, we show the cases \eqref{estimate3} holds. It suffices to show
$$ \sum_{\Nx\sim\Nd\sim N}\sum_{\Lx\sim N^4\Nn} \frac{N^{3/2}}{\langle N_1 \rangle^{3/4}L_1^{1/2+\epsilon}L_2^{1/2-\epsilon}}\Ln^{1/2}N^{-2}\Ld^{1/2} \lesssim 1 $$
Performing the $L$ summations, we reduce to
$$ \sum_{\Nx\sim\Nd\sim N} \frac{N^{-1/2}}{\langle N_1 \rangle^{3/4}} (N^4\Nn)^\epsilon \lesssim 1 $$
which is easily verified with about $N^{-1/2}$ to spare. This completes the proof for \eqref{bilinear1}.\\
The proof of \eqref{bilinear2} and \eqref{bilinear3} are very similar to the preceding one. In general, the same computation shows
$$ \|uv\|_{L^2(\R\times\R)} \lesssim \|u\|_{X^{-\alpha,1/2-\epsilon}_{\tau=\xi^5}}\|v\|_{X^{\beta,1/2+\epsilon}_{\tau=\xi^5}}. $$
for $ \alpha <2 $ and $ \alpha - \beta \leq 3/4$. We omit the detail here.
\end{proof}

\subsection*{Proof of the trilinear estimate}
In order to reduce the trilinear estimate we use the following lemma.
\begin{lemma}\label{compositionTT*}\text{[Tao \cite{Tao2001}, Lemma 3.7 Composition and TT*]}
If $k_1,k_2 \geq 1$, and $m_1, m_2$ are functions on $\R^{k_1}$ and $\R^{k_2}$ respectively, then
\begin{align}\label{composition}
\|m_1&(\xi_1,\cdots,\xi_{k_1})m_2(\xi_{k_1+1},\cdots,\xi_{k_1+k_2})\|_{[k_1+k_2;\R]} \\
  & \leq \|m_1(\xi_1,\cdots,\xi_{k_1})\|_{[k_1+1;\R]}\|m_2(\xi_{1},\cdots,\xi_{k_2})\|_{[k_2+1;\R]}. \nonumber
\end{align}
As a special case we have the $TT^*$ identity
\begin{equation}\label{TT*}
\|m(\xi_1,\cdots,\xi_{k})\overline{m(-\xi_{k+1},\cdots,-\xi_{2k})}\|_{[2k;\R]} =
   \|m(\xi_1,\cdots,\xi_{k})\|^2_{[k+1;\R]}
\end{equation}
for all functions $m:\R^k \rightarrow \R$.
\end{lemma}
For simplicity we prove the most interesting case $s=3/4$. For the first term it suffices to show that
\begin{equation*}
\left\|\frac{(\xi_1+\xi_2+\xi_3)^3\langle \xi_4\rangle^{3/4}}{\langle\tau_4-\xi_4^5\rangle^{1-b}
\prod_{j=1}^3\langle\xi_j\rangle ^s\langle\tau_j-\xi_j^5\rangle^b}\right\|_{[4,\,{\R}\times{\R}]}\lesssim 1.
\end{equation*}
Estimating $|\xi_1+\xi_2+\xi_3|$ by $\langle \xi_4 \rangle$, and
$$\langle\xi_4\rangle^{3/4+3}\lesssim \langle\xi_4\rangle^{3/2}\sum_{j=1}^3\langle\xi_j\rangle^{3/4+3/2}.$$
By symmetry we reduce to
\begin{eqnarray*}
\left\|\frac{\langle\xi_1\rangle^{-3/4}\langle\xi_3\rangle^{-3/4}\langle\xi_2\rangle^{3/2}\langle\xi_4\rangle^{3/2}}{\langle\tau_4-\xi_4^5\rangle^{1-b}
\prod_{j=1}^3\langle\tau_j-\xi_j^5\rangle^{b}}\right\|_{[4,\,{\R}\times{\R}]} \lesssim 1.
\end{eqnarray*}
We may replace $\langle\tau_2-\xi_2^5\rangle^{b}$ by $<\tau_2-\xi_2^5>^{1-b}$.
By $TT^*$ identity \eqref{TT*}, the estimate is reduced to the
bilinear estimate \eqref{bilinear1}.\\
The proof of the second term \eqref{trilinear1} is very similar to the first one but we use the composition rule \eqref{composition} instead of the $TT^*$ identity. We estimate $$ \xi_1^3 \leq \xi_1^{9/4}\left(\langle\xi_2\rangle^{3/4}+ \langle\xi_3\rangle^{3/4}+\langle\xi_4\rangle^{3/4} \right)\langle\xi_4\rangle^{3/4}.$$ The third term is the same as above and so by symmetry we reduce to
$$
\left\|\frac{\langle\xi_1\rangle^{3/2}\langle\xi_3\rangle^{-3/4}\langle\xi_2\rangle^{0}\langle\xi_4\rangle^{3/4}}{\langle\tau_4-\xi_4^5\rangle^{1-b}
\prod_{j=1}^3\langle\tau_j-\xi_j^5\rangle^{b}}\right\|_{[4,\,{\R}\times{\R}]} \lesssim 1.
$$
This is verified by \eqref{bilinear1} and \eqref{bilinear2}, as well as the composition rule \eqref{composition}.\\
The fourth term in \eqref{trilinear1} is proved in the same way. Estimating
$$ \langle\xi_4\rangle^{3/4} \leq \langle\xi_4\rangle^{1/2}\big(\langle\xi_1\rangle^{1/4}+\langle\xi_2\rangle^{1/4}+\langle\xi_3\rangle^{1/4} \big), $$
and by symmetry we reduce to
$$ \left\|\frac{\langle\xi_1\rangle^{1/2}\langle\xi_2\rangle^{1/4}\langle\xi_3\rangle^{1/4}\langle\xi_4\rangle^{1/2}}{\langle\tau_4-\xi_4^5\rangle^{1-b}
\prod_{j=1}^3\langle\tau_j-\xi_j^5\rangle^{b}}\right\|_{[4,\,{\R}\times{\R}]} \lesssim 1. $$
This is verified by \eqref{bilinear3} and $TT^*$ identity \eqref{TT*} after minorizing one of $b$ by $1-b$. \\
Finally, the third term in \eqref{trilinear1} automatically follows since it is a linear combination of other three. This conclude the proof of Proposition~\ref{trilinear}.
\begin{remark}
The trilinear estimate \eqref{trilinear1} fails for $s< \frac{3}{4}$. The counter example introduced by Kenig, Ponce and Vega \cite{KPV96} in the context of the modified KdV equation extends to here. In the frequency space, set
$$ A =\{ (\tau,\xi) \in \R^2| N \leq \xi \leq N+ N^{-3/2}, |\tau-\xi^5|\leq 1\}, $$
and $$ -A= \{(\tau,\xi) \in \R^2 | -(\tau,\xi) \in A\}. $$
Defining $\widetilde{f}(\tau,\xi) = \chi_A +\chi_{-A}$, we obtain
$$ |\widetilde{f}*\widetilde{f}*\widetilde{f}(\tau,\xi)| \gtrsim  N^{-3}\chi_R(\tau,\xi), $$
where $R$ is a rectangle located at $(N,N^5)$ of dimension $N^{-4}\times N^{5/2}$ with its longest side pointing $(1,5N^4)$ like $A$. Thus,
$$ \|\px^3(f\cdot f\cdot f) \|_{X^{s,b-1}} \gtrsim N^{s-3/4} $$
and $$ \|f\|_{X^{s,b}} \lesssim  N^{s-3/4}, $$
then \eqref{trilinear1} implies $s \geq 3/4$. This example holds good for other nonlinear terms in \eqref{trilinear1}.
\end{remark}
\begin{remark}\label{low order term}
In our general equation \eqref{fifthm} we omitted the term $u^4\px u$ from $\eqref{fifthmkdv}$. Since the term $u^4\px u$ is a lower order term, it is easier to handle than other third order terms. Once we have the 5-linear estimate 
\begin{equation}\label{5-linear}
 \|\px(u_1u_2u_3u_4u_5)\|_{X^{s,b-1}} \lesssim \prod_{j=1}^5 \|u_i\|_{X^{s,b}}, 
\end{equation}
we can insert it into the iteration. The proof of \eqref{5-linear} is similar to the preceding one. Using Lemma~\ref{compositionTT*} we reduce to two trilinear estimates and each trilinear estimate is again reduced to two bilinear estimates. The resulting bilinear estimates are supposedly easier than those in Proposition~\ref{bilinear} since there are fewer derivatives and more $u$'s. In fact, it is true for $s$ lower than $\frac{3}{4} $. 
\end{remark}

\section{Ill-posedness }\label{sec4}
In this section we give the proof of Theorem~\ref{illposed}. For simplicity, we pretend the nonlinear term is
$$ F(u) = \px^3(u^3). $$
The general case $F(u) = c_1\px^3(u^3) + c_2u^2\px^3 + c_3u\px u\px^2$ ($c_j$'s are real numbers) follows in the same manner. Our method is to approximate the fifth mKdV solution by the cubic NLS solution. This is originally introduced by Christ, Colliander and Tao \cite{CCT2003} for the mKdV equation. This method extends to the fifth order equation without substantial change. \\
Having two solutions to the cubic NLS breaking the uniform continuity of the flow map for $s<0$, we find approximate solutions to the fifth mKdV exhibiting the same property. First, we state the ill-posedness for the cubic NLS in \cite{CCT2003}.
\begin{theorem}\label{illposedNLS}
Let $  s < 0 $. The solution map of the initial value problem of the cubic NLS \eqref{cNLS} fails to be uniformly continuous. More precisely, for $ 0 < \delta \ll \epsilon \ll 1$ and $T>0$ arbitrary, there are two solutions $ u_1,u_2$ to \eqref{cNLS} satisfying \eqref{illposed1}, \eqref{illposed2} and \eqref{illposed3}.\\
Moreover, For any fixed $K \geq 1$, we can find such solutions to satisfy
\begin{equation}\label{ubound}
\sup_{0\leq t <\infty }\|u_j\|_{H^K_x} \lesssim \epsilon
\end{equation} for $j=1,2$.
\end{theorem}

\begin{remark}
Theorem~\ref{illposedNLS} is stated for the defocusing cubic NLS. The method in \cite{CCT2003} exhibiting the phase decoherence holds good for the focusing case, too. But previously another method for the focusing case was presented by Kenig, Ponce and Vega \cite{KPV2001}. They used the Galilean invariance on the soliton solutions. In our focusing case (for instance, $F(u)= -\px^3(u^3)$) one could employ their counterexample to approximate.
\end{remark}

Now we start to find the approximate solution to the fifth order mKdV equation using the NLS solutions.
Let $ u(s,y) $ solve the cubic NLS equation \eqref{cNLS}. We also assume that
$$ \sup_{0\leq t < \infty} \|u(t)\|_{H^k_x} \lesssim \epsilon $$
for a large k.
Using the change of variable
$$ (s,y) := \Big(t,\, \frac{x}{(10N^3)^{1/2}} + \sqrt{\frac{5}{2}}N^{5/2}t \Big), $$
we define the approximate solution
\begin{equation}\label{apsolution}
U_{ap}(t,x) := \frac{2}{\sqrt{3N^3}}\text{Re}\,\, e^{iNx}e^{iN^5t}u(s,y),
\end{equation}
where $N \gg 1$. \\
We want to show that $U_{ap}$ is an approximate solution to the fifth mKdV equation. A direct computation shows that
\begin{align*}
(&\partial_t -\px^5)U_{ap}(t,x) \\
= &\frac{2}{\sqrt{3N^3}} \text{Re}\,\, \Big\{e^{iNx}e^{iN^5t}\Big(\partial_su + i\py^2u + \frac{1}{\sqrt{10}N^{5/2}}\py^3u - \frac{i}{20N^5}\py^4u - \frac{1}{(10N^3)^{5/2}}\py^5u \Big)\Big\}
\end{align*}
and that
\begin{align*}
\px^3(U_{ap}^3) & = \Big(\frac{2}{\sqrt{3N^3}}\Big)^3 \frac 34 \px^3\Big\{ \text{Re}\,\, e^{iNx}e^{iN^5t}|u|^2u +\frac{1}{3}\text{Re}\,\, e^{iNx}e^{iN^5t} u^3  \Big\} \\
 & = \Big(\frac{2}{\sqrt{3N^3}}\Big)^3 \frac 34 \Bigg\{\text{Re}\,\, (iN)^3e^{iNx}e^{iN^5t}|u|^2u + \text{Re}\,\, \frac{3(iN)^2}{\sqrt{10}N^{3/2}}e^{iNx}e^{iN^5t}\py(|u|^2u) \\
 & \qquad + \text{Re}\,\, \frac{3iN}{10N^3}e^{iNx}e^{iN^5t}\py^2(|u|^2u) + \text{Re}\,\, \frac{1}{(10N^3)^{3/2}}e^{iNx}e^{iN^5t}\py^3(|u|^2u) \\
  & \qquad + \text{Re}\,\, \frac{(3iN)^3}{3}e^{3iNx}e^{3iN^5t}u^3 + \text{Re}\,\, \frac{(3iN)^2}{\sqrt{10}N^{3/2}}e^{3iNx}e^{3iN^5t}\py(u^3) \\
  & \qquad + \text{Re}\,\, \frac{3iN}{10N^3}e^{3iNx}e^{3iN^5t}\py^2(u^3) + \text{Re}\,\, \frac{1}{3(10N^3)^{3/2}}e^{3iNx}e^{3iN^5t}\py^3(u^3)  \Bigg\}.
\end{align*}
Since $u(s,y)$ is a solution of \eqref{cNLS}, three terms of the preceding equations canceled and it results in
$$ (\partial_t -\px^5)U_{ap}(t,x) + \px^3(U_{ap}^3) = E $$
where the error term $E$ is a linear combination of the real and imaginary parts of the following:
\begin{align*}
  &E_1:= N^{-4}e^{iNx}e^{iN^5t}\py(|u|^2u),\quad E_2:=N^{-11/2}e^{iNx}e^{iN^5t}\py^2(|u|^2u),\\
  &E_3:= N^{-18/2}e^{iNx}e^{iN^5t}\py^3(|u|^2u),\quad E_4:=N^{-4}e^{3iNx}e^{3iN^5t}\py(u^3), \\
  &E_5:= N^{-11/2}e^{3iNx}e^{3iN^5t}\py^2(u^3), \quad E_6:=N^{-18/2}e^{iNx}e^{iN^5t}\py^3(u^3),\\
  &E_7:= N^{-3/2}e^{3iNx}e^{3iN^5t}u^3.
\end{align*}
Next, we find the bound of the error.
\begin{lemma}\label{errorbound}
For each $j=1,\cdots,7$, let $e_j$ be the solution to the initial problem
$$ (\partial_t - \px^5)e_j = E_j; \qquad e_j(0)=0 $$
Let $\eta(t)$ be a smooth time cut-off function taking $1$ on $[0,1]$ and compactly supported. Then
$$ \|\eta(t)e_j\|_{X^{3/4,b}} \lesssim \epsilon N^{-5/2+\delta} $$ for arbitrarily small $\delta>0$.
\end{lemma}
For the proof we use the estimate of high-frequency modulations of smooth functions.
\begin{lemma}[\cite{CCT2003} Lemma 2.1]\label{modulation}
Let $-1/2 <s, \sigma \in \R^+ $ and $u\in H^\sigma(\R)$. For any $M>1, \tau\in \R^+, x_0\in\R$, and $A>0$ let
$$ v(x)= Ae^{iMx}u(\frac{x-x_0}{\tau}). $$
(i) Suppose $s\geq 0$. Then we have
$$ \|v\|_{H^1} \lesssim_s |A|\tau^{1/2}M^s\|u\|_{H^s} $$
for all $u,A,x_0$ and $M\cdot\tau \geq 1$. \\
(ii) Suppose that $s<0$ and that $\sigma \geq|s|$. Then we have
$$ \|v\|_{H^s} \lesssim_{s,\sigma} |A|\tau^{1/2}M^s\|u\|_{H^\sigma} $$
for all $u,A,x_0$ and $M^{1+(s/\sigma)}\cdot\tau \geq 1$.
\end{lemma}

\begin{proof}[Proof of Lemma~\ref{errorbound}]
Using \eqref{Xlinear} and Plancherel theorem we obtain
\begin{align*}
\|\eta(t)e_j\|_{X^{3/4,b}} & \lesssim \|\eta(t)E_j\|_{X^{3/4,b-1}} \\
  & = \|\langle\tau-\xi^5\rangle^{b-1}\langle\xi\rangle^{3/4}\widetilde{\eta(t)E_j}\|_{L^2_{\tau,\xi}} \\
  & \leq \|\langle\xi\rangle^{3/4}\widetilde{\eta(t)E_j}\|_{L^2_{\tau,\xi}}   \qquad (\because b-1 <0) \\
  & = \|\eta(t)\langle\xi\rangle^{3/4}\widehat{E_j}(t,\xi)\|_{L^2_{t,\xi}} \\
  & \leq \|\langle\xi\rangle^{3/4}\widehat{E_j}(t,\xi)\|_{L^\infty_tL^2_x([0,1]\times\R)}
\end{align*}
Thus, we reduce to show $$ \sup_{0\leq t\leq 1} \|E_j\|_{H^{3/4}_x} \lesssim \epsilon N^{-5/2+\delta}. $$
 $E_1,\cdots,E_6$ have enough negative powers of N. The above bound for these terms is obtained by \eqref{ubound}, Lemma~\ref{modulation} and the fact that $H^k$ is closed under multiplication for $k\geq 1$.
For the last term $E_7$ since there is not enough of a negative power on $N$, we need to use the fact that the modulation $e^{3iNx}e^{3iN^5t}$ is away from the the curve $\tau=\xi^5$.\\
A direct computation leads that
$$ \widetilde{\eta(t)E_7}(\tau,\xi) = N^{-3/2}\widetilde{\eta\,u^3}\Big(\tau-a,\sqrt{10}N^{3/2}(\xi-3N)\Big)\,\sqrt{10}N^{3/2} $$
where $ a= 3N^5-3\sqrt{\frac 52}N^5 +\sqrt{\frac 52}N^4\xi $. \\
Let $P_{\lambda,\mu} $ be the Littlewood-Paley projection with dyadic numbers $\lambda, \mu$.
\eqref{ubound} and the fact that $\eta(t)$ is compactly supported yield
$$ \|\widetilde{P_{\lambda,\mu}\eta\,u^3}(\tau,\xi) \|_{L^2_{\tau,\xi}} \lesssim \frac{\epsilon}{\langle\lambda\rangle^K\langle\mu\rangle^K}  $$
and so
$$ \|\widetilde{P_{\lambda,\mu}\eta\,u^3}(\tau-a,N^{3/2}(\xi-3N)) \|_{L^2_{\tau,\xi}} \lesssim N^{-3/4} \frac{\epsilon}{\langle\lambda-a\rangle^K\langle\mu-3N\rangle^K}. $$
Rewriting $\|\eta(t)E_7\|_{X^{3/4,b-1}} $ by dyadic decompositions,
\begin{align*}
&\|\eta(t)E_7\|^2_{X^{3/4,b-1}}  \\
&\lesssim \sum_{\substack{\lambda,\mu \geq 1\\ dyadic}} \langle\lambda-\mu^5\rangle^{2(b-1)}\langle\mu\rangle^{3/2} N^{-3}\Big\|\widetilde{P_{\lambda,\mu}\eta\,u^3}\Big(\tau-a,\sqrt{10}N^{3/2}(\xi-3N)\Big)\,\sqrt{10}N^{3/2}\Big\|^2_{L^2_{\tau,\xi}} \\
&\lesssim \sum_{\substack{\lambda,\mu \geq 1\\ dyadic}} \langle\lambda-\mu^5\rangle^{2(b-1)}\langle\mu\rangle^{3/2} N^{-3/2}\frac{\epsilon^2}{\langle\lambda- a\rangle^{2K}\langle\mu-3N\rangle^{2K}} \\
&\lesssim \epsilon^2 N^{10(b-1)}
\end{align*} by choosing $K$ large enough. We used the fact that $e^{3iNx}e^{3iN^5t}$ is away from the curve in the frequency space at the last inequality. Therefore, choosing $ b>\frac 12$ sufficiently close to $\frac 12$ we conclude
$$ \|\eta(t)E_7\|^2_{X^{3/4,b-1}} \lesssim \epsilon N^{-5/2+\delta}. $$

\end{proof}
Finally, we state the following perturbation result from the local well-posedness.
\begin{lemma}\label{error}
Let $u$ be a Schwartz solution to the fifth order modified KdV equation \eqref{fifthm} and v be a Schwartz solution to the approximate fifth mKdV equation
$$ \partial_t v -\px^5v + \px^3(v^3) = E $$
for some error function $E$. Let $e$ be the solution to the inhomogeneous problem
$$ \partial_t e -\px^5 e = E,\qquad e(0)=0. $$
Suppose that $$ \|u(0)\|_{H^{3/4}_x}, \|v(0)\|_{H^{3/4}_x} \lesssim \epsilon; \qquad \|\eta(t)e\|_{X^{3/4,b}} \lesssim \epsilon $$
Then we have
$$ \|\eta(t)(u-v)\|_{X^{3/4,b}} \lesssim \|u(0)-v(0)\|_{H^{3/4}} + \|\eta(t)e\|_{X^{3/4,b}} $$
In particular, we have
$$ \sup_{0\leq t\leq 1} \|u(t)-v(t)\|_{H^{3/4}} \lesssim \|u(0)-v(0) \|_{H^{3/4}} + \|\eta(t)e\|_{X^{3/4,b}}. $$
\end{lemma}

\begin{proof}
The proof is very similar to that of Lemma 5.1 in \cite{CCT2003}. Here, we give only a sketch. Writing the integral equation for $v$ with a time cut-off function $\eta(t) $
$$ \eta(t)v(t) = \eta(t)e^{t\px^5}v(0) -\eta(t)e(t) +\eta(t)\int_0^t e^{(t-t')\px^5}\px^3(v^3)(t')dt'. $$
we use \eqref{Xlinear}, \eqref{trilinear1} and a continuity argument (assuming $\epsilon $ is sufficiently small) to obtain
$$ \|\eta(t)v\|_{X^{3/4,b}} \lesssim \epsilon. $$
We repeat the same argument on the difference of the two $ w= u-v $ to get the desired result.
\end{proof}

\begin{proof}[Proof of Theorem~\ref{illposed}]
Let $0<\delta \ll \epsilon \ll 1$ and $T>0$ be given. From Theorem~\ref{illposedNLS} we can find two global solutions $u_1,u_2$ satisfying
\begin{gather}
\label{illposedNLS1}\|u_j(0)\|_{H^s_x} \lesssim \epsilon \\
\label{illposedNLS2}\|u_1(0) -u_2(0)\|_{H^s_x} \lesssim \delta \\
\label{illposedNLS3}\sup_{0\le t\le T} \|u_1(t)-u_2(t)\|_{H^s_x} \gtrsim \epsilon \\
\label{ubound1}\sup_{0\leq t \leq \infty} \|u_j(t)\|_{H^k_x} \lesssim \epsilon
\end{gather} for $s<0$ and $k\geq 6$ to be chosen later.
Define $U_{ap,1}$ and $U_{ap,2}$ as in \eqref{apsolution}, and let $U_1,U_2$ be smooth global solutions with initial data $U_{ap,1},U_{ap,2}$, respectively. Now we rescale these solutions to make them satisfy \eqref{illposed1}, \eqref{illposed2}, \eqref{illposed3}. Set
$$ U^\lambda_j(t,x) := \lambda U_j(\lambda^5t,\lambda x) $$
and similarly,
$$ U_{ap,j}^\lambda(t,x) := \lambda U_{ap,j}(\lambda^5t,\lambda x), $$
for $j=1,2$.
Then $$ U^\lambda_j(0,x) = \lambda \frac{2}{\sqrt{3N^3}}\text{Re}\,\, e^{iN\lambda x}u(0, \lambda x /(10N^3)^{1/2}). $$
From Lemma~\ref{errorbound} and Lemma~\ref{error} we have
$$ \sup_{0\leq t\leq 1} \|U_1(t)-U_2(t)\|_{H^{3/4}_x} \lesssim \epsilon N^{-5/2+\delta}. $$
An induction argument on time interval up to $\log N$ yields
\begin{equation}\label{error1}
\sup_{0\leq t \lesssim_\eta \log N} \|U_1(t)-U_2(t)\|_{H^{3/4}_x} \lesssim \epsilon N^{-5/2+\eta}
\end{equation}
for any $\eta > \delta >0$.
Applying Lemma~\ref{modulation} when $ s\geq 0$ we obtain
$$ \|U^\lambda_j(0)\|_{H^s_x} \lesssim \lambda^{s+1/2}N^{s-3/4}\|u_j(0)\|_{H^s_x}, $$
while for $s<0$, we use Lemma~\ref{modulation} (ii) for sufficiently large $k$ to obtain
$$ \|U^\lambda_j(0)\|_{H^s_x} \lesssim \lambda^{s+1/2}N^{s-3/4}\|u_j(0)\|_{H^k_x}. $$
Setting $$ \lambda:= N^{\frac{3/4-s}{1/2+s}}, $$ and from \eqref{illposedNLS1}, \eqref{ubound1}
we have \eqref{illposed1} for $U_j^\lambda(0)$. Similarly, we also get \eqref{illposed2} for $U_1^\lambda(0)-U_2^\lambda(0)$ from \eqref{illposedNLS2}. \\
Next, we show  \eqref{illposed3}. From \eqref{illposedNLS3} one can find $0<t_0 $ such that
$$ \|u_1(t_0)-u_2(t_0)\|_{L^2_x} \gtrsim \epsilon. $$ Using Lemma~\ref{modulation} we obtain
$$ \|U_{ap,1}(t_0/\lambda^5) -U_{ap,2}(t_0/\lambda^5) \|_{H^s_x} \gtrsim \lambda^{1/2+s}N^{s-3/4} \epsilon \sim \epsilon. $$
On the other hand, using the hypothesis $s>-\frac{7}{24}$ and \eqref{error1}
$$  \|U^\lambda_{ap,j}(t)-U^\lambda_j(t)\|_{H^s} \lesssim \lambda^{\max(0,s) +1/2}\epsilon N^{-5/2+\eta} \lesssim \epsilon $$
for $ 0<t\lesssim_\eta \log N/\lambda^5$ and sufficiently small $\eta>0$. A triangle inequality shows
$$ \|U^\lambda_1(t_0/\lambda^5) - U^\lambda_2(t_0/\lambda^5) \|_{H^s_x} \gtrsim\epsilon $$
for $t_0/\lambda^5 \ll \log N/\lambda^5$. Choosing $\lambda$(and hence N) large enough that $t_0/\lambda^5 <T$, we get \eqref{illposed3}.
This completes the proof.

\end{proof}

\end{document}